\documentclass[12pt]{amsart}
\usepackage{latexsym, amsmath, amssymb, amsthm, amsfonts, amstext, amscd, amsxtra, amsbsy, mathrsfs, stmaryrd, mathdots, bbold}
\usepackage[all]{xy}
\usepackage{tikz}
\usepackage{soul}
\usepackage[makeroom]{cancel}
\usepackage{bbm}
\usepackage{color}
\usepackage{hyperref}
\usepackage{shorttoc}

\newtheorem{theorem}{Theorem}[section]
\newtheorem{lemma}[theorem]{Lemma}
\newtheorem{prop}[theorem]{Proposition}

\theoremstyle{definition}

\newtheorem{example}[theorem]{Example}
\newtheorem{remark}[theorem]{Remark}

\newcommand{\mbf}{\mathbf}
\newcommand{\mbb}{\mathbb}

\def\mf{\mathfrak}

\title[Application of Schur-Weyl duality to Springer theory]{Application of Schur-Weyl duality to Springer theory}

\author[ Z. Dong and H. Ma]{Zhijie Dong and Haitao Ma}
\address{(H.Ma )College of mathematics science, Harbin Engineering University, Harbin, 150001, China.}
\address{(Z.Dong )Harbin Institute of Technology, Harbin, 150001, China.}
\email{dongmouren@gmail.com(Z. Dong)\\
 hmamath@hrbeu.edu.cn(H. Ma)
 }

\date{\today}

\begin{document}

\begin{abstract}
In \cite{FMX19}, it is proved that the convolution algebra of  top  Borel-Moore homology on Steinberg variety of type $B/C$ realizes  $U(sl_n^{\theta})$, where $sl_{n}^{\theta}$ is the fixed point subalgebra of involution on $sl_n$. So  top Borel-Moore homology of the partial Springer's fibers gives the representations of $U(sl_n^{\theta})$. In this paper, we study these representations  using the Schur-Weyl duality and Springer theory.
\end{abstract}

\maketitle

\section{Introduction}

In \cite{G91}, V. Ginzburg observed that one can realize a irreducible representation of $U(sl_n)$ using the top Borel-Moore homology of partial springer fibers of type A. A. Braverman and D. Gaitsgory  gave another interpretation of the action  using the  Schur-Weyl duality between $U(sl_n)$ and Weyl group $S_d$ and  Springer theory \cite{BG99}.

Let $\mf g$ be a Lie algebra and $\theta: \mf g\rightarrow \mf g$ be an involution on $\mf g$.
Let $\mf g^{\theta}$ be the fixed point subalgebra of $\mf g$.
In \cite{SS99}, M. Sakamoto and T. Shoji gave the Schur-Weyl duality of  Ariki-Koike algebras which includes $sl_n^{\theta}$ as a special case.
In \cite{BKLW14} , the authors considered the convolution algebra on double flag varieties associated to the algebraic group $O_{2r+1}(\mbb F_q)$
( $Sp_{2r}(\mbb F_q)$),
which is $q$-Schur type algebra.
The geometric realizations of the corresponding $\imath$-quantum group (quantization of $U(\mathfrak g^{\theta}$))are hence obtained. They also gave the Schur-Weyl duality between the $\imath$-quantum group and Hecke algebra of type $B/C$. In \cite{FMX19}, it was proved that the convolution algebra of the top  Borel-Moore homology on Steinberg variety of type $B/C$ realizes  $U(sl_n^{\theta})$ , so the Borel-Moore homology of partial Springer's fibers give a representation of $U(sl_n^{\theta})$.

This paper is devoted to study the representation of $U(sl_n^{\theta})$  given by the Borel-Moore homology of partial Springer's fibers of type C. The type B case can be given similarly so we omit it. More precisely, let $a \in \mathcal O_A$ be the nilpotent orbit corresponding to a Young diagram $A$ of type C.  We prove in
theorem \ref{Main theorem}
top homology of partial Springer fibers over $a$
can be naturally identified with a representation of $U(sl_n^{\theta})$, which is determined by Schur-Weyl duality and Springer correspondence.

This paper is organized as follows. In section 2 and section 3, we recall the Schur-Weyl duality of  Ariki-Koike algebras and the Springer theory. In section 4, we give the proof of  the main result of this paper and an example to explain it.

\section{Schur-Weyl duality}
Let    $g=gl_{m_1} \oplus gl_{m_2}$ be a Levi subalgebra of $gl_m$ with $ m = m_1 + m_2$.
Let $W_{d}$ be  the complex reflection group $S_d \ltimes (\mathbb{Z}/2\mathbb{Z})^d$.
The group $W_{d}$ is generated by $s_1,s_2,\cdots,s_d$,
where $s_2,\cdots,s_d$ are the generators of $S_d$ corresponding to transpositions $(1,2),\cdots,(n-1,n)$,
and $s_1$ satisfies the relation $s_1^2 = 1,(s_1s_2)^2 = (s_2s_1)^2$,$s_1s_i = s_is_1$ for $i \geq 3$.
 Let $V = \mathbb{C}^m$ be the natural representation of $U(gl_m)$.
We can define  a left $U(g)$ action and right $W_{d}$  action  on $V^{\otimes d}$ as follows.
The action of $U(g)$ on $V^{\otimes d}$ is obtained by the restriction of $U(gl_m)$.
Let $\varepsilon = \{v_1^{1},\cdots v_{m_1}^1, v_{1}^{2}\cdots,v_{m_2}^2\}$ be the natural basis of $V$.
 Define a function $b : \varepsilon \rightarrow \mathbb{N} $ by $b(v_j^i) = i$.
 The group $S_d$ acts on $V^{\otimes d}$ by permuting the components of the tensor product, while $s_1$ acts on $V^{\otimes n}$ by
 $$s_{1}(x_1\otimes \cdots \otimes x_n) = (-1)^{b(x_1)}x_1\otimes \cdots \otimes x_n,$$
 where $x_i \in \varepsilon$.
 Then Schur-Weyl duality holds between the left $U(g)$ action and the right $W_{d}$  action on $V^{\otimes d}$.

 Since $U(g) = U(gl_{m_1})\otimes  U(gl_{m_2})$, irreducible representations of $U(g)$ are parameterized by 2-tuples $\mathbf{\lambda} = (\lambda^1, \lambda^2)$ of  Young diagrams $\lambda^i$ with $l(\lambda^i) \leq m_i$,
 where $l(\lambda^i)$ is the number of rows of $\lambda^i$.
  Let $\Lambda_{m_1,m_2}$ be the set of 2-tuple $\mathbf{\lambda}$ of  Young diagrams such that  $l(\lambda^i) \leq m_i$, $|\lambda^1| + |\lambda^2| = d$, where  $|\lambda^i|$ is the number of boxes in  $\lambda^i$.
     Let $V_{\mathbf{\lambda}}$   be the irreducible $U(g)$-module corresponding to $\lambda$.
   We recall that the  irreducible representation of $W_d$ are in one to one correspondence to the ordered pairs of Young diagrams $(\mu, \nu)$ where $\mu = (\mu_i), \mu_1 \geq \mu_2 \geq \cdots \geq \mu_r$ is a partition of $d_1$ and $\nu =(\nu_j),\nu_1 \geq \nu_2 \geq \cdots \geq \nu_s$ is a partition of $d_2$ such that $d_1+d_2 =d$.
   Let $Z_{(\mu,\nu)}$   be the irreducible $W_{d}$-module corresponding to $(\mu,\nu)$.

  \begin{prop}\cite{SS99}\label{Schur-Weyl}
 The $ U(g)\otimes W_{d}$-module $V^{\otimes d}$ decomposes as
  $$V^{\otimes d} = \bigoplus_{\mathbf{\lambda} \in \Lambda_{m_1,m_2}}Z_{\mathbf{\lambda}} \otimes V_{\mathbf \lambda}.$$
\end{prop}

\section{Springer theory}

Let $N = 2n+1$, $D =2d$, $N>D$.
 Let $V $ be $\mbb C^{2d}$  with a nondegenerate skew-symmetric bilinear form $( , )$, and $G = Sp(V)$.
 Let
$$\mathcal{Q}_{N, D}=\{ \mathbf{d}=(d_i) \in \mbb N^{N} \mid   d_i = d_{N+1-i},\quad \textstyle \sum_{i=1}^{N} d_i = D \}.$$
For any $U\subseteq V$, let $U^{\perp} = \{x \in V| (x,y)=0, \ \forall y\in U\}$.
For any $\mbf d \in \mathcal{Q}_{N, D}$, we  set partial flag
$$\mathcal{F}_{\mathbf{d}}=\{ F=( 0=V_0 \subset V_1 \subset \cdots \subset V_{N}=V)\ \mid \  V_i=V_{N-i}^{\perp} ,\ \text{dim}( V_{i}/V_{i-1})= d_i,\  \forall i \}. $$
Denote
$$\mathcal{C}=\{ F=( 0=V_0 \subset V_1 \subset \cdots \subset V_{D}=V)\ \mid \  V_i=V_{N-i}^{\perp} ,\ \text{dim}( V_{i}/V_{i-1})= 1,\  \forall i \}$$
for convenience.
Let $\mathcal{F} = \sqcup _{\mathbf{d} \in \mathcal{Q}_{N,D}}\mathcal{F}_{\mathbf{d}}$.

The nilpotent cone of $\mathfrak{g} = \mathfrak{sp}_{D}$, denote by $\mathcal{N}$, consists of all  nilpotent element in $\mathfrak g$ .
Let
$\mathcal N_{N,D} = \{a \in \mathfrak{g} | a^n = 0\}$,
and $\widetilde{\mathcal N_{N,D}} = T^*\mathcal F $ (resp. $\widetilde{\mathcal N} = T^*\mathcal C $ ) be the cotangent bundle of $\mathcal F$(resp. $\mathcal C$).
 More precisely, $\widetilde{\mathcal N_{N,D}}$ can be identified with the set of all pairs
\begin{equation*}
\{(F,x) \in \mathcal F \times \mathfrak{g}\mid x(F_i) \subseteq F_{i-1}, ~ \forall i \},
\end{equation*}
and $\widetilde{\mathcal N}$ can be identified with the set of all pairs
\begin{equation*}
\{(F,x) \in \mathcal C \times \mathfrak{g} \mid x(F_i) \subseteq F_{i-1}, ~ \forall i \}.
\end{equation*}

For $\mathbf{d} \in \mathcal{Q}_{N,D}$,  denote  by $\widetilde{\mathcal N_{N,D}}^{\mathbf{d}}$  the $\mathbf d$'s component of $\widetilde{\mathcal N_{N,D}}$.
Let $\pi^{N,D}: \widetilde{\mathcal N_{N,D}} \rightarrow \mathfrak{g}, (F,x) \mapsto x$ be projection map.
 Let $\mu: \widetilde{\mathcal N} \rightarrow \mathcal N, (F,x) \mapsto x$ be the springer resolution of  $\mathcal N$.  Let $\tilde{\mathfrak{g}}$ be the Grothendieck simultaneous resolution
 \begin{equation*}
\{(F,x) \in \mathcal C \times \mathfrak{sp}_{2d} \mid x(F_i) \subseteq F_{i}, ~ \forall i \}.
\end{equation*}
 Let $\mu_{\tilde{\mathfrak{g}}}: \tilde{\mathfrak{g}} \rightarrow \mathfrak{sp}_{2d}, (F,x) \mapsto x$ be  projection map.
 Let $\mathfrak{g}_{rs}$ be  regular semisimple part of $\mathfrak{g}$.
  We have the following cartesian diagram
$$
\CD
  \tilde{\mathfrak{g}}_{rs} @> >> \tilde{\mathfrak{g}} @< << \tilde{\mathcal{N}} \\
  @V \mu_{rs} VV @V \mu_{\tilde{\mathfrak{g}}} VV @V \mu VV  \\
  \mathfrak{g}_{rs} @> >> \mathfrak{g} @< << \mathcal{N}
\endCD
$$
 The springer sheaf for $ G $, denoted by $Spr$, is the perverse sheaf
$$Spr = \mu_{*} \underline{\mathbb{C}}_{\widetilde{\mathcal N}}[dim \widetilde{\mathcal N}].$$
 Then there is an isomorphism
$$\mathbb{C}[W] \simeq \mathrm{End}(Spr),$$
where $W = W_d$.

 For any $x_0 \in \mathfrak{g}_{rs}$, there is  a surjective map $\pi_{1}(\mathfrak{g}_{rs}, x_0)\twoheadrightarrow W$.
   Let
 $$L_{rs}: \mathbb{C}[W]-\mathrm{mod} \rightarrow \mathrm{Loc} (\mathfrak{g}_{rs},\mathbb{C})$$
 be the map that assigns  $V$ to the local system corresponding to it.
 For any  $\mathbb{C}[W]$-module $\rho$,
 let  $\mathcal{S}_{\rho}$ be   IC sheaf $\mathrm{IC}(\mathfrak{g}_{rs},L_{rs}(\rho))$ on $\mathfrak{g}$.
 Let us choose a $G-$equivariant isomorphism between  vector spaces $\mathfrak g$ and $ \mathfrak g^{*}$, then one can regard the Fourier-Laumon transform as  a functor
 $$Four_{\mathfrak g} : D_{\mathbb{C}^* \times G}^b(\mathfrak g, \mathbb{C})\rightarrow D_{\mathbb{C}^* \times G}^b(\mathfrak g, \mathbb{C}).$$

 Abuse of notation, for any irreducible $W$-representation $\rho = (\lambda,\mu)$,
let $\nu = (\nu_i)$ be the sequence defined by $\nu_{2i-1} = \mu_i (1 \leq i \leq s)$, $\nu_{2i} = \lambda_i (1 \leq i \leq r)$,
and $\nu_j = 0$ otherwise.  Define Young diagram $A_{\rho} = (a_i)$ of type C by the following equalities,

(i) if $\nu_i \geq \nu_{i+1}$, then $a_i = 2\nu_i$,

(ii) if $\nu_i = \nu_{i+1} - 1$, then $a_i = 2\nu_i+1, a_{i+1} = 2\nu_i +1$,

(iii) if $\nu_i \leq \nu_{i+1} - 2$, then $a_{i} = 2\nu_{i+1}-2, a_{i+1} = 2\nu_{i} +2$.
By  \cite[Theorem 3.3]{S79}, the springer correspondence can be given  by the following map
$$\rho = (\lambda ,\mu) \rightarrow (A_{\rho},\phi_{\rho}),$$
where $\phi_{\rho}$  denote the local system on the nilpotent orbit correspondence to $\rho$ defined explicitly in  \cite{S79}.
\begin{remark}
Note that in section 2, we already have a correspondence between bipartitions of $n$ and irreducible representations of $W$. This two correspondences are equal up to duality of Young diagrams.
\end{remark}

We have the following well-known theorem.
\begin{theorem}
Let $h : \mathcal{N} \hookrightarrow \mathfrak{g}$ be the inclusion map. For any irreducible  representation $\rho$ of $\mathbb{C}[W]$, we have the following isomorphism of the perverse sheaves
$$Four_{\mathfrak g}(\mathcal{S}_{\rho}) \cong \mathrm{IC}_{A_{\rho},\phi_{\rho}}.$$
\end{theorem}

\section{Main results}

 We define a perverse sheaf $\mathcal{L}$ such that for any $\mathbf{d} \in \mathcal{Q}_{N,D}$ $$\mathcal{L}_{\mathbf d} = \mathcal{L}|_{\widetilde{\mathcal N_{N,D}}^{\mathbf d}}= \underline{\mathbb C}[\mathrm{dim}(\widetilde{\mathcal N_{N,D}}^{\mathbf d})].$$
 Let $\widetilde{\mathfrak g_{N,D}}$  be  the variety of all pairs
\begin{equation*}
\{(F,x) \in \mathcal F \times \mathfrak{g}\mid x(F_i) \subseteq F_{i}, ~ \forall i \}.
\end{equation*}
Define $p : \widetilde{\mathfrak g_{N,D}} \rightarrow \mathfrak g,(F,x) \mapsto x$.
Let $\widetilde{\mathfrak g_{N,D}}^{\mathbf d}$ be the $\mathbf d$'s component of $\widetilde{\mathfrak g_{N,D}}$.
We have
$\mathrm{dim} \widetilde{\mathfrak g_{N,D}}^{\mathbf d} = \mathrm{dim} (\mathfrak{g}).$
 Define
$$\mathcal{K}= \underline{\mathbb C}_{\widetilde{\mathfrak g_{N,D}}}[\mathrm{dim}(\mathfrak{g})].$$

\begin{lemma}
There is a canonical  isomorphism
$$Four_{\mathfrak g}(\pi^{N,D}(\mathcal{L})) \simeq p_* \mathcal{K}.$$
\end{lemma}
\begin{proof}
 Oberseve $\widetilde{\mathcal{N}_{N,D}}$ and  $\widetilde{\mathfrak{g}}$ are the subbundles of the trivial bundle $\mathcal{F} \times \mathfrak{sp}_{2d}$. Under the $G$ -quivariant  $\mathfrak{sp}_{2d}^{*}\simeq \mathfrak{sp}_{2d}$.
We have
  $$\widetilde{\mathcal{N}_{N,D}} =  \widetilde{\mathfrak{sp}_{2d}}^\bot.$$
  The lemma follows from the \cite[Corollary 6.8.11]{Ac}.
\end{proof}

Let $E$ be $\mathbb{C}^{N}$. Recall that $E^{\otimes d}$ is a $W$-module, and
 $\mathcal{S}_{E^{\otimes d}} = \mathrm{IC}(\mathfrak{g}_{rs},L_{rs}(E^{\otimes d})).$
 \begin{prop}\cite[Proposition 2.7.10]{Ac} \label{Acbookprop}
 Let $X$ and $Y$ be connected, locally path-connected, and semilocally simply connected topological spaces, and let $f:(X,x_0) \rightarrow (Y,y_0)$ be a covering map. For $\mathfrak{F} \in Loc(X,\mathbb K)$, there is a natural isomorphism
 $$Mon_{y_0}(^{\circ}f_{!}\mathfrak F ) \cong \mathbb{K}[\pi_1(Y,y_0)] \otimes_{\mathbb{K}[\pi_1(X,x_0)]}Mon_{x_0}(\mathfrak F).$$
 \end{prop}

\begin{lemma} \label{lemma2}
There is a canonical  isomorphism
$$p_*{\mathcal K} \simeq \mathcal S_{E^{\otimes d}}.$$

\end{lemma}
\begin{proof}
Since the map $p$ is small, $p_{*}\mathcal K$ is equal to the intersection cohomology extension of its restriction local system on $\mathfrak{g}_{rs}$.
We only need to prove
$$p_*{\mathcal K}|_{\mathfrak g_{rs}} \simeq \mathcal S_{E^{\otimes d}} |_{\mathfrak g_{rs}}.$$
Denote $p_{rs} = p|_{p^{-1}(\mathfrak g_{rs})}$.
For any $\mathbf{d} =(d_i) \in \mathcal{Q}_{N,D}$,
let $W_{\mathbf d}$ be the group $ S_{d_1} \times S_{d_2} \times \cdots \times S_{d_{n}} \times (S_{d_n+1} \ltimes (\mathbb{Z}/2\mathbb{Z})^{d_{n+1}})$.
The map $p_{rs}$ is a covering map with Galois group $W/W_{\mathbf d}$.
Denote by $\mathbb{C}[X]$ be $\mathbb{C}$-vector space generated by the set $X$.
By proposition \ref{Acbookprop}, $p_*{\mathcal K}|_{\mathfrak g_{rs}}$ is a local system corresponding to  a natural $W$-module on $\mathbb{C}[W/W_{\mathbf d}]$.
More precisely,
fix $t \in \mathfrak{g}_{rs}$.
There exist a basis $\{e_1,e_2,\cdots,e_D\}$ such that $\mathfrak{g}_{rs}$ is the diagonal matrix corresponding to this basis.
 Denote $X_{\mathbf{d}}= p^{-1}(t) \cap \widetilde{\mathfrak g_{N,D}}^{\mathbf d}$.  Define a set
$$\Theta_{\mathbf d}=\{A= (a_{ij}) \in \mathrm{Mat}_{N \times D}(\mathbb{N}) | \sum_{i,j} a_{ij} = 2d,a_{ij} = a_{N+1-i,N +1 -j}, \sum_{i}a_{ij} = 1,\sum_{j}a_{ij} = d_i\}.$$
For any $A \in \Theta_{\mathbf d}$, we can construct an element
$$ 0 \subset V_1 \subset V_2 \subset \cdots \subset V$$
 in $X_{\mathbf{d}}$,
  where $V_i = \sum_{j\in[1,N],a_{ij} = 1} \mathbb{C}e_j$.
  Then we have $X_{\mathbf d}$ and $\Theta_{\mathbf d}$ are in bijection.
 The set $\Theta_{\mathbf d}$ has a $W$-action defined as follows.
 For any $\sigma \in W$,
    $$\sigma((V_i)) =(\sigma(V_i)),$$
   where $ \sigma(V_i)= \sum_{j\in[1,N],a_{ij} = 1} \mathbb{C}e_{\sigma(j)}.$
    Then $\mathbb{C}[X_{\mathbf d}]$  is  a W-module, and $\mathbb{C}[\Theta_{\mathbf d}]$ is also a W-module induced by the bijection between the two sets.
 Let
  $$F = 0 \subset U_1 \subset U_2 \subset \cdots \subset V \in X_{\mathbf d},$$
where $U_i = \sum_{j=1}^{\sum_{m= 1}^id_i} \mathbb{C}e_j$.
Define
   $$\varphi : \mathbb{C}[W/W_{\mathbf d}] \rightarrow \mathbb{C}[X_d], w \mapsto w(F).$$
We can check $\varphi$ is a W-module isomorphism.

Next consider $p^{-1}(t) = \bigcup_{\mathbf{d} \in \mathcal{Q}_{N,D}}X_{\mathbf{d}}$ which consists of all components. It is in bijection to the union of $\Theta_d$ which is denoted by
$$\Theta =\{A= (a_{ij}) \in \mathrm{Mat}_{N \times N}(\mathbb{D}) | \sum_{i,j} a_{ij} = 2d,a_{ij} = a_{N+1-i,N +1 -j}, \sum_{i}a_{ij} = 1 \}.$$
 For any $A \in  \Theta$, every column have only one non-zero number 1, and the first d columns can decide the last d columns.  So the number of the element in $P^{-1}(t)$ is $N^d$.
 The dimension of the $E^{\otimes d}$  is also $N^d$.
 Define
 $$\chi: \mathbb{C}[\Theta] \rightarrow E^{\otimes d},A =(a_{ij}) \mapsto f_{i_1} \otimes f_{i_2} \otimes \cdots \otimes f_{i_d},$$
  where $i_k$ is the number such that $a_{k,i_k} = 1$,$\{f_1, \cdots, f_N\} $ is the natural basis of $E$.
 It is easy to see $\chi$ is injective and $W$-equivariant. So $\chi$ is an isomorphism as $W$-module since $dim \mathbb{C}[\Theta] = dim E^{\otimes d} = N^d$.
Recall $p_{*}(\mathcal K)|_{\mathfrak g_{rs}} = L_{rs}(\mathbb{C}[\Theta])|_{\mathfrak g_{rs}}$ and $ S_{E^{\otimes d}}|_t = L_{rs}(E^{\otimes d})$.
 Then we have
$$p_{*}(\mathcal K)|_{\mathfrak g_{rs}} \simeq \mathcal S_{E^{\otimes d}}|_{\mathfrak g_{rs}}.$$
 The lemma follows.
\end{proof}

Let $A$ be the Young diagram of type C,  and $\mathcal{O}_{ A}$ be the nilpotent orbit corresponding to $A$. Let $\mathcal{N}_{N,D}^{\mathbf d}$ be  image of $\pi^{N,D}$ on $\mathbf{d}$'s component.
For  $a \in \mathcal{O}_{ A}$,
we don't know the dimension of $(\pi^{N,D})^{-1}(a)\cap \mathcal{F}_{\mathbf{d}}$.
 By definition of semismall map,
 the dimension is lower than $codim_{\mathcal{N}_{N,D}^{\mathbf d}}(\mathcal O_{A})$.
So we call  $codim_{\mathcal{N}_{N,D}^{\mathbf d}}(\mathcal O_{A})$-dimensional homology the top homology.
We have
$$H_{top}((\pi^{N,D})^{-1}(a) )= \bigoplus_{\mathbf{d} \in \mathcal{Q}_{N,D}} H_{codim_{\mathcal{N}_{N,D}^{\mathbf d}}(\mathcal O_{A})}((\pi^{N,D})^{-1}(a)\cap \mathcal{F}_{\mathbf{d}}),$$
where $H_{*}$ is the Borel-Moore homology.

For  $\rho =(\lambda,\mu)$, denote $\check{\rho} = (\check{\lambda},\check{\mu})$, where $\check{\bullet}$ is the dual of Young diagram $\bullet$.
 The following theorem is the main theorem of this paper.
\begin{theorem} \label{Main theorem}
For any Young diagram $A$ of type C, $a \in \mathcal{O}_{A}$, then the space $H_{top}((\pi^{N,D})^{-1}(a))$ can be naturally identified with the space of the representation $\oplus _{\{\rho\in Irr(W)|A_{\rho} = A\}} V_{\check{\rho}}$,
where the summation is over irreducible representations $\rho$ such that $A_{\rho} =A$.
\end{theorem}

\begin{proof}
Let $\pi^{N,D}_{\mathbf d}$ be the restriction of $\pi^{N,D}$ on $\mathbf d$'s component.
 Since $\pi^{N,D}|_{\mathbf d}$ is a semismall \cite{BM83},
so  the pushforward of  a perverse sheaf is also a perverse sheaf.
 By the decomposition theorem,
 there is a canonical  isomorphism
  $$\pi^{N,D}_{\mathbf d, *}(\mathcal{L}_{\mathbf d})
  \simeq \bigoplus_{(B,\iota)} IC_{(B,\iota)} \otimes W(B,\iota)_{\mathbf d},$$
  where $\iota$ is the local system on $\mathcal{O}_B$.
By proper base change,  $(\pi^{N,D}_{\mathbf d, *} \mathcal L_{\mathbf d})_a \cong R\Gamma (\pi^{N,D}_{\mathbf d}(a), \underline{\mathbb{C}})[dim \mathcal N_{N,D}^{\mathbf d}]$.
We  obtain
$$H^i(\pi^{N,D}_{\mathbf d}(a),\underline{\mathbb{C}}) \cong  \bigoplus_{(B,\iota)} H^{i- dim \mathcal N_{N,D}^{\mathbf d}}(IC_{B,\iota}|_{\mathcal O_B})_a \otimes W(B,\iota)_{\mathbf d}. $$
Let $i = codim_{\mathcal{N}_{N,D}^{\mathbf d}}(\mathcal O_{A})$. Then
$$H^{codim_{\mathcal{N}_{N,D}^{\mathbf d}}(\mathcal O_{A})}(\pi^{N,D}_{\mathbf d}(a),\underline{\mathbb{C}}) \cong  \bigoplus_{(B,\iota)} H^{-dim \mathcal O_1}(IC_{B,\iota}|_{\mathcal O_A})_a \otimes W(B,\iota)_{\mathbf d}.$$
By the property of IC, we have
               $$  H^{-dim \mathcal O_A}(IC_{(B,\iota),a}) =\left\{\begin{array}{ll}
\mathbb{C}& \text{if} \ A =B; \\[.15in]
0& \text{otherwise}.
 \end{array}   \right.
$$
By the definition of the Borel-Moore homology,
$$H^{codim_{\mathcal{N}_{N,D}^{\mathbf d}}(\mathcal O_{A})}((\pi^{N,D}_{\mathbf d})^{-1}(a),\underline{\mathbb{C}}) \cong H_{codim_{\mathcal{N}_{N,D}^{\mathbf d}}(\mathcal O_{A})}((\pi^{N,D}_{\mathbf d})^{-1}(a)).$$
Then
$$ H_{codim_{\mathcal{N}_{N,D}^{\mathbf d}}(\mathcal O_{A})}((\pi^{N,D})^{-1}(a)\cap \mathcal{F}_{\mathbf{d}}) \cong \bigoplus_{\iota} W(A,\iota)_{\mathbf d}.$$
So
$$ H_{top}((\pi^{N,D})^{-1}(a)) \cong \bigoplus_{\mathbf d \in \in \mathcal{Q}_{N,D}}\bigoplus_{\iota} W(A,\iota)_{\mathbf d}.$$
Combine all components together, we have
$$\pi^{N,D}_{ *}(\mathcal{L})
  \simeq \bigoplus_{(B,\iota)} IC_{(B,\iota)} \otimes W(B,\iota).$$
So
$$W(B,\iota) = \bigoplus_{\mathbf d \in \in \mathcal{Q}_{N,D}} W(B,\iota)_{\mathbf d},$$
and
$$H_{top}((\pi^{N,D})^{-1}(a)) \cong \bigoplus_{\iota} W(A,\iota).$$

Since $Four_{\mathfrak g}$ is an involution,
we have the following canonical isomorphism
 $$\pi^{N,D}_{*}(\mathcal L) \simeq Four_{\mathfrak g}(\mathcal S_{E^{\otimes d}}).$$
 By the proposition \ref{Schur-Weyl}, we have
 $$\mathcal S_{E^{\otimes d}} = \bigoplus_{\rho \in irr(W)} \mathcal S_{\rho} \otimes V_{\rho}.$$
Then
  $$\bigoplus_{(A,\iota)} IC_{(A,\iota)} \otimes W(A,\iota) \simeq \pi^{N,D}_{*}(\mathcal{L}) \simeq \bigoplus_{\rho \in irr(W)} Four_{\mathfrak g}(\mathcal S_{\rho}) \otimes V_{\rho} =  \bigoplus_{\rho \in irr(W)} IC_{A_{\rho},\phi_{\rho}} \otimes V_{\rho}.$$
 So we have
               $$ W(A,\iota) =\left\{\begin{array}{ll}
V_{\check{\rho}}& \text{if} \ (A,\iota) = (A_{\rho},\phi_{\rho}); \\[.15in]
0& \text{otherwise}.
 \end{array}   \right.
$$
The theorem follows.
\end{proof}

\begin{theorem}
The action of $sl_{n+1}\bigoplus gl_n$ on
$V_{\rho}$ via Schur-Weyl duality coincides with the  action by convolution.
\end{theorem}

\begin{proof}
Denote $Z$ the Steinberg variety $\widetilde{\mathcal N_{N,D}}\times_{\mathfrak{g}}\widetilde{\mathcal N_{N,D}}$.
We define an action of $H(Z)$ on the sheaf $\mathfrak{F}= p_*(\underline{\mathbb C}_{\widetilde{\mathfrak g_{N,D}}}[\mathrm{dim}(\mathfrak{g})])$ geometrically.
This means that for any open set $U\subset \mathfrak{g}$, define an action
of $H(Z)$  on $\mathfrak{F}(U)=H(\pi^{-1}(U))$. This is done by convolution.
In the paper\cite{FMX19}, it is proved that there is a map $U(sl_N^{\theta}) \rightarrow H(Z)$,
hence we have $sl_N^{\theta}$ acting on the sheaf $\mathfrak{F}$, which we denoted by $G$.
We next define another action of $sl_N^{\theta}$ acting on the sheaf $\mathfrak{F}$.
Recall in lemma \ref{lemma2},
since the map $p$ is small, the sheaf $\mathfrak{F}$ is the IC sheaf associated with the local system $E^{\otimes d}$. Using the isomorphism between $sl_N^{\theta}$ and $sl_{n+1}\bigoplus gl_{n+1}$,
we have the natural action of
$sl_N^{\theta}$ on $E^{\otimes d}$.
By IC continuation principle, this action, denoted by $A$, extends to the sheaf $\mathfrak{F}$, and to show the action $A$ coincides with the action $G$, we only need to show they coincide on the regular semisimple part $\mathfrak{g}_{rs}$.
For a  regular semisimple element $x\in \mathfrak{g}_{rs}$, the fiber $p^{-1}(x)$.
is analysed in lemma \ref{lemma2} and by result in \cite[Section 6]{BKLW14}, letting q=1 yields our claim that $A$ coincides with $G$.
Now we decompose the sheaf $\mathfrak{F}= \bigoplus_{\rho \in Irr(W) }  \mathcal S_{\rho}\otimes V_{\rho}$.
Algebraically, the space $V_{\rho}$ gets the action of $sl_N^{\theta}$ by Schur-Weyl duality and it is the action of  $sl_N^{\theta}$ on $E^{\otimes d}$ restricted on the subspace $V_{\rho}$.
Geometrically, the space $V_{\rho}$ is identified with the top degree cohomology of the stalk at $x\in \mathcal{N}$ and by definition of the action $G$, it is given by convolution.
Since we proved that $A$ coincides with $G$, the action of $sl^{\theta}$ on
$V_{\rho}$ via Schur-Weyl duality coincides with the action by convolution.
\end{proof}

\begin{remark}
The action of $sl_N^{\theta}$ on $E^{\otimes d}$ defined in \cite[Section 6]{BKLW14} looks slightly different from the action of $sl_{n+1}\oplus gl_{n}$ in  \cite{SS99} after identifying $sl_N^{\theta}$ and $sl_{n+1}\oplus gl_{n}$, but they are the same after a change of basis.
\end{remark}

\begin{example}
Consider the case $n = 2, d = 2$.
 Let $\theta : sl_5 \rightarrow sl_5$ be the involution defined by $\theta(e_i) =f_{5- i},\theta(f_i) = e_{5-i}$. As we know, $sl_5^{\theta}$ is generated by
 $$E_1 = e_1 + f_4,
 E_2 = e_2 + f_3,
 F_1 = f_1 + e_4,
 F_2 = f_2 + e_3,
H_1 = h_1 - h_4 ,
 H_2 = h_2 - h_3. $$
Let $V$ be the 5-dimensional $\mathbb{C}$-vector space with basis $\{e_1,\cdots, e_5\}$.
The Schur-Weyl duality between $sl_5^{\theta}$ and $\mathbb{C}[S_2\ltimes \mathbb{Z}_2^2]$  is given as follows.
The action of $\mathbb{C}[S_2\ltimes \mathbb{Z}_2^2]$ on $V^{\otimes 2}$ is given by
  $$s_1(e_i\otimes e_j) = e_j \otimes e_i,[1]_1(e_i \otimes e_j)=e_i \otimes e_{6-j},$$
 where $s_1$ is the generator of $S_2$,
 $[1]_1 = ([1],0)$ in  $\mathbb{Z}_2^2$.
   There are four 1-dimensional  modules and one  2-dimensional simple module.
   The Young diagram of type C  consists of four cases $(1,1,1,1),(2,1,1),(2,2),$
   $(4)$.
 By the Springer correspondence of type C, we have:\\
\begin{center}
\begin{tabular}{|c|c|c|c|}
\hline
Irr(W) &dim& ($\lambda,\mu$) & Young diagram\\
\hline
 Sign & 1&((1,1),0) & (1,1,1,1) \\
\hline
 Ssign & 1&(0,(1,1)) & (2,1,1)  \\
\hline
 Lsign & 1&(2,0) & (2,2)\\
\hline
reguler & 2&(1,1) & (2,2)\\
\hline
triv & 1&(0,2) & (4)  \\
\hline
\end{tabular}
\end{center}
The partition of type C in this case is
  \begin{align*}
 \mathcal{Q}_{5, 4}=& \{\mathbf d_1=(1,1,0,1,1),\mathbf d_2=(0,1,2,1,0),\mathbf d_3=(1,0,2,0,1),\\
     &  \mathbf d_4 = (0,2,0,2,0),\mathbf d_5 =(2,0,0,0,2),\mathbf d_6=(0,0,4,0,0)\}.
\end{align*}
For any $x \in \mathcal O_A$,
by direct computation,
the following table gives the fiber of $x$ in different component.
 \begin{center}
\begin{tabular}{|c|c|c|c|c|c|c|}
\hline
   &$\mathbf d_1$& $\mathbf d_2$ & $\mathbf d_3$ & $\mathbf d_4$ & $\mathbf d_5$ & $\mathbf d_6$\\
\hline
 (4) & pt&$\emptyset$ &$\emptyset$ & $\emptyset$&$\emptyset$ & $\emptyset$ \\
\hline
 (2,2) & $\mathbb{P}_1 \cup \mathbb{P}_1\cup \mathbb{P}_1$&$pt_1 \cup pt_2$ &$pt_1 \cup pt_2$  & pt &pt &$\emptyset$  \\
\hline
(2,1,1)& X  & pt&pt & $\mathbb{P}_1$& $\mathbb{P}_1$& $\emptyset$\\
\hline
(1,1,1,1) & $\mathcal F_{\mathbf d_1}$ & $\mathcal F_{\mathbf d_2}$ & $\mathcal F_{\mathbf d_3}$ & $\mathcal F_{\mathbf d_4}$ & $\mathcal F_{\mathbf d_5}$ & pt\\
\hline
\end{tabular}
\end{center}
where $X$  is an irreducible variety of dimension 2.

If $x \in \mathcal O_{(2,1,1)}$, the following table gives dimension of top homology  in different component.
 \begin{center}
\begin{tabular}{|c|c|c|c|c|c|c|}
\hline
   &$\mathbf d_1$& $\mathbf d_2$ & $\mathbf d_3$ & $\mathbf d_4$ & $\mathbf d_5$ & $\mathbf d_6$\\
\hline
(2,1,1)& X  & pt&pt & $\mathbb{P}_1$& $\mathbb{P}_1$& $\emptyset$\\
\hline
 $dim H_{top}$& 1 & 0 & 0 & 1& 1 & 0\\
\hline
\end{tabular}
\end{center}
In the $\mathbf d_i$'s component, where $i = 2,3$. Since $codim_{N^{\mathbf d_i}_{5,4}}(\mathcal{O}_{x}) = \frac{6-4}{2} = 1$,
we have
$H_{top}(\pi_{\mathbf d_i}^{-1}(x)) =H_1(pt) =0$.
 So we have $dim H_{top}(\pi^{-1}(x)) = 3$.
In the other cases, the top dimension we defined is  really the top dimension of the fiber.

\end{example}

\end{document}